\documentclass[12pt]{article}

\usepackage[margin=1.5cm]{geometry}

\usepackage{cite}

\usepackage{tikz}

\usepackage{mathtools}
\usepackage{graphicx}
\usepackage{amssymb,amsmath,wasysym}
\usepackage{amsfonts}
\usepackage{comment}
\usepackage{enumerate}

\newtheorem{Remark}{Remark}

\newcommand{\bLambda}{\mbox{\boldmath$\Lambda$}}

\newcommand{\subR}{\mbox{$\scriptstyle\mathbb{R}$}}

\begin{document}
	

\title{Random walk on a quadrant: mapping to a one-dimensional level-dependent Quasi-Birth-and-Death process (LD-QBD)}

\author{
%
%
Ma{\l}gorzata M. O'Reilly
\thanks{Ma{\l}gorzata M. O'Reilly is supported by the Australian Research Council Discovery Project DP180100352.}
~~\thanks{Discipline of Mathematics, University of Tasmania, Tas 7001, Australia, email: malgorzata.oreilly@utas.edu.au}
\and 
Zbigniew Palmowski
\thanks{Zbigniew Palmowski was partially supported by the National Science Centre (Poland) under the grant 2021/41/B/HS4/00599.}
~~\thanks{Faculty of Pure and Applied Mathematics, Wroc{\l}aw University of Science and Technology, 50-370 Wroc{\l}aw, Poland, email: zbigniew.palmowski@pwr.edu.pl}
\and
Anna Aksamit
\thanks{Anna Aksamit is supported by the Australian Research Council Early Career Researcher Award DE200100896.}
~~\thanks{School of Mathematics and Statistics, The University of Sydney, NSW 2006, Australia, email: anna.aksamit@sydney.edu.au}
}

\date{\today}

\maketitle

\section{Introduction}
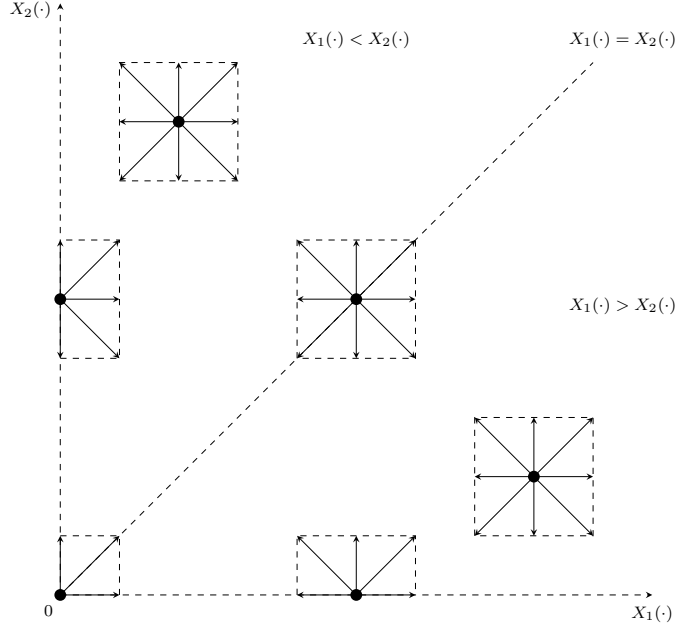
\begin{figure}[h!]
	\begin{center}
		\resizebox{0.5\textwidth}{!}{
		\begin{tikzpicture}[>=stealth,redarr/.style={->}]

		\draw [ dashed, ->] (0,0) -- (10,0);
		\draw [ dashed, ->] (0,0) -- (0,10);
		\draw [dashed] (0,0) -- (9,9);
		
		\draw (-0.5,10) node[anchor=north, below=-0.17cm] {\scriptsize{\color{black} $X_2(\cdot)$}};
		
		\draw (10,-0.2) node[anchor=north, below=-0.17cm] {\scriptsize{\color{black} $X_1(\cdot)$}};
		
		\draw (-0.2,-0.2) node[anchor=north, below=-0.17cm] {\scriptsize{\color{black} $0$}};
		
		\draw (9.5,9.5) node[anchor=north, below=-0.17cm] {\scriptsize{\color{black} $X_1(\cdot)=X_2(\cdot)$}};
		
		\draw (5,9.5) node[anchor=north, below=-0.17cm] {\scriptsize{\color{black} $X_1(\cdot)<X_2(\cdot)$}};
		
		\draw (9.5,5) node[anchor=north, below=-0.17cm] {\scriptsize{\color{black} $X_1(\cdot)>X_2(\cdot)$}};

		\node at (8,2) [black,circle,fill,inner sep=2pt]{};
		\draw [black, ->] (8,2) -- (8,3);
		\draw [black, ->] (8,2) -- (8,1);
		\draw [black, ->] (8,2) -- (9,2);
		\draw [black, ->] (8,2) -- (7,2);
		\draw [black, ->] (8,2) -- (9,3);
		\draw [black, ->] (8,2) -- (7,1);
		\draw [black, ->] (8,2) -- (7,3);
		\draw [black, ->] (8,2) -- (9,1);
		\draw [dashed] (7,1) -- (7,3) -- (9,3) -- (9,1) -- (7,1);

		\node at (2,8) [black,circle,fill,inner sep=2pt]{};
		\draw [black, ->] (2,8) -- (3,8);
		\draw [black, ->] (2,8) -- (1,8);
		\draw [black, ->] (2,8) -- (2,9);
		\draw [black, ->] (2,8) -- (2,7);
		\draw [black, ->] (2,8) -- (3,9);
		\draw [black, ->] (2,8) -- (1,7);
		\draw [black, ->] (2,8) -- (3,7);
		\draw [black, ->] (2,8) -- (1,9);
		\draw [dashed] (1,7) -- (3,7) -- (3,9) -- (1,9) -- (1,7);

		\node at (5,5) [black,circle,fill,inner sep=2pt]{};
		\draw [black, ->] (5,5) -- (6,5);
		\draw [black, ->] (5,5) -- (4,5);
		\draw [black, ->] (5,5) -- (5,6);
		\draw [black, ->] (5,5) -- (5,4);
		\draw [black, ->] (5,5) -- (6,6);
		\draw [black, ->] (5,5) -- (4,4);
		\draw [black, ->] (5,5) -- (6,4);
		\draw [black, ->] (5,5) -- (4,6);
		\draw [dashed] (4,4) -- (4,6) -- (6,6) -- (6,4) -- (4,4);

		\node at (5,0) [black,circle,fill,inner sep=2pt]{};
		\draw [black, ->] (5,0) -- (6,0);
		\draw [black, ->] (5,0) -- (4,0);
		\draw [black, ->] (5,0) -- (6,1);
		\draw [black, ->] (5,0) -- (4,1);
		\draw [black, ->] (5,0) -- (5,1);
		\draw [dashed] (4,0) -- (4,1) -- (6,1) -- (6,0);

		\node at (0,5) [black,circle,fill,inner sep=2pt]{};
		\draw [black, ->] (0,5) -- (0,6);
		\draw [black, ->] (0,5) -- (0,4);
		\draw [black, ->] (0,5) -- (1,6);
		\draw [black, ->] (0,5) -- (1,4);
		\draw [black, ->] (0,5) -- (1,5);
		\draw [dashed] (0,4) -- (1,4) -- (1,6) -- (0,6);
		
		\node at (0,0) [black,circle,fill,inner sep=2pt]{};
		\draw [black, ->] (0,0) -- (0,1);
		\draw [black, ->] (0,0) -- (1,0);
		\draw [black, ->] (0,0) -- (1,1);
		\draw [dashed] (0,1) -- (1,1) -- (1,0);

		\end{tikzpicture}
	}
	\end{center}
	\caption{Neighbourhood random walk on a quadrant $\{(X_1(t),X_2(t),\varphi(t)):t\geq 0\}$. Given current state $(n,m,i)$, at the next jump in the process we may observe: (i) a change in the phase $\varphi(\cdot)$ only without a change in the position $(X_1(\cdot),X_2(\cdot))$, with $(n,m,i)\to(n,m,j)$, $j\not= i$; or (ii) a change in both, with $(n,m,i)\to (n',m',i)$, $(n,m)\not= (n',m')$, $j\not= i$. The distribution of time spent in the position $(n,m)$ is phase-type with parameters that depend on $(n,m)$.
}
	\label{fig:RWquadrant}
\end{figure}

We consider a neighbourhood random walk on a quadrant, $\{(X_1(t),X_2(t),\varphi(t)):t\geq 0\}$, with state space
\begin{eqnarray*}
	\mathcal{S}&=&\{(n,m,i):n,m=0,1,2,\ldots;i=1,2,\ldots,k(n,m)\}
\end{eqnarray*}
forming a single communicating class, see Figure~\ref{fig:RWquadrant}. We refer to $X_1(t)$ and $X_2(t)$ as the position variables and to $\varphi(t)$ as the environment variable. The process evolves in time as follows.

	Assuming start in state $(n,m,i)$, the process spends exponentially distributed amount of time in $(n,m,i)$ according to some parameter $\lambda_i^{(n,m)}$.
	

	Upon leaving state $(n,m,i)$ the process moves to some state $(n^{'},m^{'},j)$ with $j\in\{1,\ldots,k(n^{'},m^{'})\}$ and $n^{'}\in\{n-1,n,n+1\}$, $m^{'}\in\{m-1,m,m+1\}$, according to some probabilities $(p_{n;a}^{m;b})_{i,j}$ with $a,b\in\{+,-,0\}$.
%
	 


We transform this process into a one-dimensional LD-QBD $\{(Z(t),\chi(t)):t\geq 0\}$ with level variable $Z(t)$ and phase variable $\chi(t)$. Using this transform we find its transient and stationary analysis using matrix-analytic methods in Ramaswami~\cite{Rama}, Joyner and Fralix~\cite{2016JF}, and Phung-Duc et al.~\cite{phung2010simple}, as well as the distribution at first hitting times.

We generalise a spatially-coherent uniformisation of an SFM to a QBD proposed in~\cite{2013BO_uni} and construct a sequence of two-dimensional LD-QBDs that converge in distribution to a two-dimensional stochastic fluid model $\{(Y_1(t),Y_2(t),\varphi(t)):t\geq 0\}$, which describes a movement on a quadrant in which the position changes in a continuous manner according to rates $dY_1(t)/dt=c_{1,\varphi(t)}$ and $dY_2(t)/dt=c_{2,\varphi(t)}$. 

\newpage
\section{Random walk as a 2-D LD-QBD}

We start by formally defining the two-dimensional random walk $\{(X_1(t),X_2(t),\varphi(t)):t\geq 0\}$
that we consider in this paper. We assume that the random walk evolves in time under the following assumptions.
\begin{itemize}
	\item Assuming $(X_1(t),X_2(t),\varphi(t))=(n,m,i)$ for some $n,m\in\{0,1,2,\ldots\}$ and $i\in\mathcal{S}_{n,m}=\{1,\ldots,k(n,m)\}$ at some time $t\geq 0$, the process spends exponentially distributed amount of time $\tau(i)$ in $i$ according to distribution $\tau(i)\sim Exp(\lambda_i^{(n,m)})$, with parameters $\lambda_i^{(n,m)}$ collected in a matrix ${\bLambda}^{(n,m)}=diag(\lambda_i^{(n,m)})_{i\in\{1,\ldots,k(n,m)\}}$.
	\item Upon leaving phase $i$ the process moves to some phase $j\not= i$, $j\in\{1,\ldots,k(n^{'},m^{'})\}$ and simultaneously moves to some position $(X_1(t+\tau(i)),X_2(t+\tau(i)))=(n^{'},m^{'})$, $n^{'}\in\{n-1,n,n+1\}$, $m^{'}\in\{m-1,m,m+1\}$, according to probabilities $(p_{n;a}^{m;b})_{i,j}$, $a,b\in\{+,-,0\}$ collected in a matrix ${\bf P}_{n;a}^{m;b}=[(p_{n;a}^{m;b})_{i,j}]_{i\in\{1,\ldots,k(n,m)\},j\in\{1,\ldots,k(n^{'},m^{'})\}}$, such that
	\begin{eqnarray*}
		&&m^{'}=m+1 \mbox{ when }b=+,\\
		&&m^{'}=m \mbox{ when }b=0,\\
		&&m^{'}=[m-1]^+ \mbox{ when }b=-,\\
		&&n^{'}=n+1 \mbox{ when }a=+,\\
		&&n^{'}=n \mbox{ when }a=0,\\
		&&n^{'}=[n-1]^+ \mbox{ when }b=-,
	\end{eqnarray*}
	where $[u]^+=max\{0,u\}$, and
	\begin{eqnarray*}
		(p_{n;0}^{m;0})_{i,i}&=&0,\\
		\sum_{a\in\{+,-,0\}}\sum_{b\in\{+,-,0\}}
		\sum_{j\in\{1,\ldots,k(n^{'},m^{'})\}} (p_{n;a}^{m;b})_{i,j}&=&1.
	\end{eqnarray*}
	\item The process $\{(X_1(t),X_2(t),\varphi(t)):t\geq 0\}$ is irreducible, and so its state space
	\begin{eqnarray*}
		\mathcal{S}&=&\{(n,m,i):n,m=0,1,\ldots;i=1,\ldots,k(n,m)\}
	\end{eqnarray*}
	forms a single communicating class.
\end{itemize}

In other words, our random walk
has transition rates 
\begin{eqnarray*}
	(\lambda_{n;a}^{m;b})_{i,j}&=&\lambda_i^{(n,m)}\times (p_{n;a}^{m;b})_{i,j}
	= [{\bLambda}^{(n,m)} \times {\bf P}_{n;a}^{m;b}]_{ij},
\end{eqnarray*}
that depend on the current state $(X_1(t),X_2(t),\varphi(t))=(n,m,i)$ at time $t$. Note that at the moment of jump in the phase process, each level variable may move one level up, stay at the same level, or move one level down (if the current level is positive). Therefore this is a neighbourhood random walk.

Denote
\begin{eqnarray*}
	\bLambda^{(n,m)}&=&diag(\lambda_i^{(n,m)})_{i\in\{1,\ldots,k(n,m)\}},\\
	\bLambda^{(n,\bullet)}&=&diag(\bLambda^{(n,m)})_{m=0,1,\ldots ,n-1},\\
	\bLambda^{(\bullet,m)}&=&diag(\bLambda^{(n,m)})_{n=0,1,\ldots ,m-1}.
\end{eqnarray*}

The process $\{(X_1(t),X_2(t),\varphi(t)):t\geq 0\}$ forms a two-dimensional level-dependent Quasi-Birth-and-Death process (LD-QBD) with two-dimensional level variable $(X_1(t),X_2(t))$ taking values in $\{0,1,2,\ldots\}\times \{0,1,2,\ldots\}$ and one-dimensional phase variable $\varphi(t)$ taking values in $\cup_{n,m=0,1,2,\ldots}\mathcal{S}_{n,m}$.

Its generator
is given by ${\bf Q} =
[{\bf Q}^{[n,m;n^{'},m^{'}]}]_{n,m,n^{'},m^{'}\in\{0,1,2,\ldots\}}$ with the block matrices 
\begin{eqnarray*}
	{\bf Q}^{[n,m;n^{'},m^{'}]}=[q_{(n,m,i)(n^{'},m^{'},j)}]_{
		i\in\{1,2,\ldots \widehat k(n,m)\},
		j\in\{1,2,\ldots \widehat k(n^{'},m^{'})\}} 
\end{eqnarray*}
such that
\begin{eqnarray*}
	q_{(n,m,i)(n^{'},m^{'},j)}&=&
	(\lambda_{n;a}^{m;b})_{i,j}1\{n^{'}=[n+a]^+,m^{'}=[m+b]^+\}
\end{eqnarray*}
for $(n,m,i)\not=(n^{'},m^{'},j)$, and $$q_{(n,m,i)(n,m,i)}=-\sum_{(n,m,i)\not=(n^{'},m^{'},j)}q_{(n,m,i)(n^{'},m^{'},j)},$$ 
where $1\{\cdot\}$ denotes an indicator function.


\section{Mapping to a 1-D LD-QBD}

One of the key observation of this paper is that one can transform above random walk into
a corresponding one-dimensional LD-QBD $\{(Z(t),\chi(t)):t\geq 0\}$ with level variable $Z(t)\geq 0$ and phase variable $\chi(t)=(\epsilon_1(t),\epsilon_2(t),\varphi(t))$. We achieve this by using the following procedure:
\begin{itemize}
	\item Let $Z(t)=\max\{X_1(t),X_2(t)\}$, and so $Z(t)=n$ when $n=\max\{X_1(t),X_2(t)\}$, for $n\in\{0,1,2,\ldots\}$.
	\item $\epsilon_1(t)\in\{0,1,2\}$, and $\epsilon_2(t)\in\{0,1,\ldots,n\}$ when $Z(t)=n$, are such that
	\begin{eqnarray}
		\epsilon_1(t)&=&
		\left\{
		\begin{array}{cc}
			0& \mbox{when }X_1(t)=X_2(t),\\
			1& \mbox{when }X_1(t)>X_2(t),\\
			2& \mbox{when }X_1(t)<X_2(t);
		\end{array}
		\right.
		\\
		\epsilon_2(t)&=& \min\{X_1(t),X_2(t)\}.
	\end{eqnarray}
\end{itemize}
Then the correspondence between $\{(X_1(t),X_2(t),\varphi(t)):t\geq 0\}$ and $\{(Z(t),\chi(t)):t\geq 0\}$ is one-to-one since the state of one process corresponds to exactly one state of another in this pair. That is, $\{(Z(t),\chi(t)):t\geq 0\}$ is an alternative way of representing $\{(X_1(t),X_2(t),\varphi(t)):t\geq 0\}$.

We have $i\in\{1,\ldots ,\widehat k(n,\epsilon_1,\epsilon_2) \}$ whenever $Z(t)=n$ and $\chi(t)=(\epsilon_1,\epsilon_2,i)$, with
\begin{itemize}
\item $\widehat k(n,\epsilon_1,\epsilon_2)=k(n,n)$ if $\epsilon_1=0$, and then $\epsilon_2\in \mathcal{E}(n,\epsilon_1)=\{n\}$,

\item $\widehat k(n,\epsilon_1,\epsilon_2)=k(n,\epsilon_2)$ if $\epsilon_1=1$, and then $\epsilon_2\in \mathcal{E}(n,\epsilon_1)=\{0,1,\ldots,n-1\}$,

\item $\widehat k(n,\epsilon_1,\epsilon_2)=k(\epsilon_2,n)$ if $\epsilon_1=2$, and then $\epsilon_2\in \mathcal{E}(n,\epsilon_1)=\{0,1,\ldots,n-1\}$.
\end{itemize}

We note that the generator ${\bf Q} =
[{\bf Q}^{[n,n^{'}]}]_{n,n^{'}\in\{0,1,2,\ldots\}}$ of the process $\{(Z(t),\chi(t)):t\geq 0\}$,
is given by
\begin{eqnarray*}
	{\bf Q} 
&=&
	\begin{bmatrix}
		{\bf Q}^{[0,0]} & {\bf Q}^{[0,1]} & {\bf O} & {\bf O} & \cdots\\
		{\bf Q}^{[1,0]} & {\bf Q}^{[1,1]} & {\bf Q}^{[1,2]} & {\bf O} & \cdots\\
		{\bf O} & {\bf Q}^{[2,1]} & {\bf Q}^{[2,2]} & {\bf Q}^{[2,3]} &  \cdots\\
		{\bf O} & {\bf O} & {\bf Q}^{[3,2]} & {\bf Q}^{[3,3]} & \cdots \\
		\vdots & \vdots & \vdots & \vdots & \ddots
	\end{bmatrix},
\end{eqnarray*}
with the block matrices ${\bf Q}^{[n,n^{'}]}$ as specified below.

For $n=0$,
\begin{eqnarray*}
	{\bf Q}^{[0,0]}&=&
	{\bf Q}^{[0,0]}_{0,0}
	=
	\bLambda^{(0,0)}{\bf P}^{0;0}_{0;0},\\
	{\bf Q}^{[0,1]}&=&
	[{\bf Q}^{[n,n^{'}]}_{0,\epsilon_1^{'}}]_{\epsilon_1^{'}\in\{0,1,2\}}
	=
	\left[
	\begin{array}{ccc}
		{\bf Q}^{[0,1]}_{0,0}& {\bf Q}^{[0,1]}_{0,1}& {\bf Q}^{[0,1]}_{0,2}
	\end{array}
	\right],
\end{eqnarray*}
with
\begin{eqnarray*}
	{\bf Q}^{[0,1]}_{0,0}&=&
	\bLambda^{(0,0)}
	{\bf P}^{0;+}_{0;+}
	,\\
	{\bf Q}^{[0,1]}_{0,1}&=&
	\bLambda^{(0,0)}
	{\bf P}^{0;0}_{0;+},\\
	{\bf Q}^{[0,1]}_{0,2}&=&
	\bLambda^{(0,0)}
	{\bf P}^{0;+}_{0;0}.
\end{eqnarray*}

For $n\geq 1$,
\begin{eqnarray*}
	{\bf Q}^{[n,n]}&=&
	[{\bf Q}^{[n,n]}_{\epsilon_1,\epsilon_1^{'}}]_{\epsilon_1,\epsilon_1^{'}\in\{0,1,2\}}
	=
	\left[
	\begin{array}{ccc}
		{\bf Q}^{[n,n]}_{0,0}&{\bf Q}^{[n,n]}_{0,1}&{\bf Q}^{[n,n]}_{0,2}\\
		{\bf Q}^{[n,n]}_{1,0}&{\bf Q}^{[n,n]}_{1,1}&{\bf O}\\
		{\bf Q}^{[n,n]}_{2,0}&{\bf O}&{\bf Q}^{[n,n]}_{2,2}
	\end{array}
	\right]
\end{eqnarray*}
with ${\bf Q}^{[n,n]}_{1,2}={\bf O}$ and ${\bf Q}^{[n,n]}_{2,1}={\bf O}$ since a change from $X_1<X_2$ to $X_1>X_2$ or from $X_1>X_2$ to $X_1<X_2$ may not occur without observing $X_1=X_2$ first, and
\begin{eqnarray*}
	{\bf Q}^{[n,n]}_{0,0}&=&
	=-\bLambda^{(n,n)}+\bLambda^{(n,n)}{\bf P}^{n;0}_{n;0},\\
	{\bf Q}^{[n,n]}_{0,1}&=&
%
\left[
\begin{array}{cccc}
	{\bf O}&\ldots&{\bf O}&
	\bLambda^{(n,n)}{\bf P}^{n;-}_{n;0}
\end{array}
\right],	
	\\
	{\bf Q}^{[n,n]}_{0,2}&=&
\left[
\begin{array}{cccc}
	{\bf O}&\ldots&{\bf O}&
	\bLambda^{(n,n)}{\bf P}^{n;0}_{n;-}
\end{array}
\right],
\end{eqnarray*}
and
\begin{eqnarray*}
	{\bf Q}^{[n,n]}_{1,0}&=&
	\left[
	\begin{array}{cccc}
		{\bf O}&\ldots&{\bf O}&\bLambda^{(n,n-1)}{\bf P}^{n-1;+}_{n;0}
	\end{array}
	\right]^T,\\
	{\bf Q}^{[n,n]}_{2,0}&=&
	\left[
	\begin{array}{cccc}
		{\bf O}&\ldots&{\bf O}&\bLambda^{(n-1,n)}{\bf P}^{n;0}_{n-1;+}
	\end{array}
	\right]^T,
\end{eqnarray*}
and
\begin{eqnarray*}
\lefteqn{
	{\bf Q}^{[n,n]}_{1,1}=
	-\bLambda^{(n,\bullet)}+\bLambda^{(n,\bullet)}
}
	\nonumber\\
	&&\times
	\begin{bmatrix}
		{\bf P}^{0;0}_{n;0} & {\bf P}^{0;+}_{n;0}& {\bf O} & \cdots & \cdots & {\bf O}\\
		{\bf P}^{1;-}_{n;0} & {\bf P}^{1;0}_{n;0} & {\bf P}^{1;+}_{n;0} & {\bf O} &\cdots & {\bf O}\\
		{\bf O} & {\bf P}^{2;-}_{n;0} & {\bf P}^{2;0}_{n;0} & {\bf P}^{2;+}_{n;0} & \ddots &  \vdots\\
		{\bf O} & \ddots & \ddots & \ddots & \ddots & {\bf O}\\
		{\bf O} & \ddots & \ddots &  {\bf P}^{n-2;-}_{n;0} & {\bf P}^{n-2;0}_{n;0} & {\bf P}^{n-2;+}_{n;0}\\
		{\bf O} & \ddots & \ddots & {\bf O} & {\bf P}^{n-1;-}_{n;0} & {\bf P}^{n-1;0}_{n;0}
	\end{bmatrix},
\end{eqnarray*}
and
\begin{eqnarray*}
	\lefteqn{
	{\bf Q}^{[n,n]}_{2,2}=
	-\bLambda^{(\bullet,n)}
	+\bLambda^{(\bullet,n)}
}
	\nonumber\\
	&&\times
	\begin{bmatrix}
		{\bf P}^{n;0}_{0;0} & {\bf P}^{n;0}_{0;+}& {\bf O} & \cdots & \cdots & {\bf O}\\
		{\bf P}^{n;0}_{1;-} & {\bf P}^{n;0}_{1;0} & {\bf P}^{n;0}_{1;+} & {\bf O} &\cdots & {\bf O}\\
		{\bf O} & {\bf P}^{n;0}_{2;-} & {\bf P}^{n;0}_{2;0} & {\bf P}^{n;0}_{2;+} & \ddots &  \vdots\\
		{\bf O} & \ddots & \ddots & \ddots & \ddots & {\bf O}\\
		{\bf O} & \ddots & \ddots &  {\bf P}^{n;0}_{n-2;-} & {\bf P}^{n;0}_{n-2;0} & {\bf P}^{n;0}_{n-2;+}\\
		{\bf O} & \ddots & \ddots & {\bf O} & {\bf P}^{n;0}_{n-1;-} & {\bf P}^{n;0}_{n-1;0}
	\end{bmatrix}.
\end{eqnarray*}

For $n\geq 1$,
\begin{eqnarray*}
	{\bf Q}^{[n,n+1]}&=&
	\left[
	\begin{array}{ccc}
		{\bf Q}^{[n,n+1]}_{;0,0}& {\bf Q}^{[n,n+1]}_{0,1}& {\bf Q}^{[n,n+1]}_{0,2}\\
		{\bf O}& {\bf Q}^{[n,n+1]}_{1,1}& {\bf O}\\
		{\bf O}& {\bf O}&{\bf Q}^{[n,n+1]}_{2,2}
	\end{array}
	\right],
\end{eqnarray*}
with
\begin{eqnarray*}
	{\bf Q}^{[n,n+1]}_{0,0}&=&
	\bLambda^{(n,n)}{\bf P}^{n;+}_{n;+},
	\\
	{\bf Q}^{[n,n+1]}_{0,1}&=&
	\left[
	\begin{array}{cccc}
		{\bf O}&\ldots&{\bf O}&\bLambda^{(n,n)}{\bf P}^{n;0}_{n;+}
	\end{array}
	\right],
	\\
	{\bf Q}^{[n,n+1]}_{0,2}&=&
	\left[
	\begin{array}{cccc}
		{\bf O}&\ldots&{\bf O}&\bLambda^{(n,n)}{\bf P}^{n;+}_{n;0}
	\end{array}
	\right],
\end{eqnarray*}
and
\begin{eqnarray*}
\lefteqn{
	{\bf Q}^{[n,n+1]}_{1,1}=
	\bLambda^{(\bullet,n)}
}
	\nonumber\\
	&&\times
	\begin{bmatrix}
		{\bf P}^{0;0}_{n;+} & {\bf P}^{0;+}_{n;+}& {\bf O} & \cdots & \cdots &\cdots & {\bf O}\\
		{\bf P}^{1;-}_{n;+} & {\bf P}^{1;0}_{n;+} & {\bf P}^{1;+}_{n;+} & {\bf O} &\cdots &\cdots & {\bf O}\\
		{\bf O} & {\bf P}^{2;-}_{n;+} & {\bf P}^{2;0}_{n;+} & {\bf P}^{2;+}_{n;+} & \ddots &  \cdots &\vdots\\
		{\bf O} & \ddots & \ddots & \ddots & \ddots &\ddots & \vdots\\
		{\bf O} & \ddots & {\bf O} & {\bf P}^{n-2;-}_{n;+} & {\bf P}^{n-2;0}_{n;+} & {\bf P}^{n-2;+}_{n;+} & {\bf O}\\
		{\bf O} & \ddots & \ddots & {\bf O} & {\bf P}^{n-1;-}_{n;+} & {\bf P}^{n-1;0}_{n;+} & {\bf P}^{n-1;+}_{n;+}
	\end{bmatrix},
	\nonumber\\
\end{eqnarray*}
and
\begin{eqnarray*}
	\lefteqn{
	{\bf Q}^{[n,n+1]}_{2,2}=
	\bLambda^{(\bullet,n)}
}
\\
&&\times
	\begin{bmatrix}
		{\bf P}^{n;+}_{0;0} & {\bf P}^{n;+}_{0;+}& {\bf O} & \cdots & \cdots &\cdots & {\bf O}\\
		{\bf P}^{n;+}_{1;-} & {\bf P}^{n;+}_{1;0} & {\bf P}^{n;+}_{1;+} & {\bf O} &\cdots &\cdots & {\bf O}\\
		{\bf O} & {\bf P}^{n;+}_{2;-} & {\bf P}^{n;+}_{2;0} & {\bf P}^{n;+}_{2;+} & {\bf O} &\cdots &  \vdots\\
		{\bf O} & \ddots & \ddots & \ddots & \ddots &\ddots & \vdots\\
		{\bf O} & \ddots & {\bf O} & {\bf P}^{n;+}_{n-2;-} & {\bf P}^{n;+}_{n-2;0} & {\bf P}^{n;+}_{n-2;+}& {\bf O}\\
		{\bf O} & \ddots & \ddots & {\bf O} & {\bf P}^{n;+}_{n-1;-} & {\bf P}^{n;+}_{n-1;0} & {\bf P}^{n;+}_{n-1;+}
	\end{bmatrix}
	.
	\nonumber\\
\end{eqnarray*}

For $n\geq 1$,
\begin{eqnarray*}
	{\bf Q}^{[n,n-1]}&=&
	\left[
	\begin{array}{ccc}
		{\bf Q}^{[n,n-1]}_{0,0}& {\bf O}& {\bf O}\\
		{\bf Q}^{[n,n-1]}_{1,0}& {\bf Q}^{[n,n-1]}_{1,1}& {\bf O}\\
		{\bf Q}^{[n,n-1]}_{2,0}& {\bf O}&{\bf Q}^{[n,n-1]}_{2,2}
	\end{array}
	\right],
\end{eqnarray*}
with
\begin{eqnarray*}
	{\bf Q}^{[n,n-1]}_{0,0}&=&
	\bLambda^{(n,n)}{\bf P}^{n;-}_{n;-} ,
	\\
	{\bf Q}^{[n,n-1]}_{1,0}&=&
	\left[
	\begin{array}{ccccc}
		{\bf O}&\ldots &{\bf O}& \bLambda^{(n,n-2)}{\bf P}^{n-2;+}_{n;-}&
		\bLambda^{(n,n-1)}{\bf P}^{n-1;0}_{n;-}
	\end{array}
	\right]^T,
	\\
	{\bf Q}^{[n,n-1]}_{2,0}&=&
	\left[
	\begin{array}{ccccc}
		{\bf O}&\ldots &{\bf O}& \bLambda^{(n-2,n)}{\bf P}^{n;-}_{n-2;+}&
		\bLambda^{(n-1,n)}{\bf P}^{n;-}_{n-1;0}
	\end{array}
	\right]^T,
\end{eqnarray*}

and
\begin{eqnarray*}
\lefteqn{
	{\bf Q}^{[n,n-1]}_{1,1}=
	\bLambda^{(n,\bullet)}
}
\\
&&\times
	\begin{bmatrix}
		{\bf P}^{0;0}_{n;-} & {\bf P}^{0;+}_{n;-}& {\bf O} & \cdots & \cdots & {\bf O}\\
		{\bf P}^{1;-}_{n;-} & {\bf P}^{1;0}_{n;-} & {\bf P}^{1;+}_{n;-} & {\bf O} &\cdots & {\bf O}\\
		{\bf O} & {\bf P}^{2;-}_{n;-} & {\bf P}^{2;0}_{n;-} & {\bf P}^{2;+}_{n;-} & \ddots &  \vdots\\
		{\bf O} & \ddots & \ddots & \ddots & \ddots & {\bf O}\\
		{\bf O} & \ddots & {\bf P}^{n-3;-}_{n;-} & {\bf P}^{n-3;0}_{n;-} & {\bf P}^{n-3;+}_{n;-} & {\bf O}\\
		{\bf O} & \ddots & \ddots & {\bf O} & {\bf P}^{n-2;-}_{n;-} & {\bf P}^{n-2;0}_{n;-}\\
		{\bf O} & \ddots & \ddots & {\bf O} & {\bf O} & {\bf P}^{n-1;-}_{n;-}
	\end{bmatrix},
\end{eqnarray*}
and
\begin{eqnarray*}
\lefteqn{
	{\bf Q}^{[n,n-1]}_{2,2}=
	\bLambda^{(\bullet,n)}
}
\\
&&\times
	\begin{bmatrix}
		{\bf P}^{n;-}_{0;0} & {\bf P}^{n;-}_{0;+}& 0 & \cdots & \cdots & {\bf O}\\
		{\bf P}^{n;-}_{1;-} & {\bf P}^{n;-}_{1;0} & {\bf P}^{n;-}_{1;+} & {\bf O} &\cdots & {\bf O}\\
		{\bf O} & {\bf P}^{n;-}_{2;-} & {\bf P}^{n;-}_{2;0} & {\bf P}^{n;-}_{2;+} & \ddots &  \vdots\\
		{\bf O} & \ddots & \ddots & \ddots & \ddots & {\bf O}\\
		{\bf O} & \ddots & {\bf P}^{n;-}_{n-3;-} & {\bf P}^{n;-}_{n-3;0} & {\bf P}^{n;-}_{n-3;+} & {\bf O}\\
		{\bf O} & \ddots & \ddots &  {\bf O} & {\bf P}^{n;-}_{n-2;-} & {\bf P}^{n;-}_{n-2;0}\\
		{\bf O} & \ddots & \ddots & {\bf O} & {\bf O} & {\bf P}^{n;-}_{n-1;-}
	\end{bmatrix}.
\end{eqnarray*}

\begin{Remark}
	Suppose that $\{(X_1(t),X_2(t),\varphi(t)):t\geq 0\}$ is a level-independent 2-D QBD , that is
	\begin{itemize}
		\item $\{1,\ldots ,k(n,m)\}=\{1,\ldots ,k\}$ for all $n,m$,
		\item $\lambda_i^{(n,m)}=\lambda_i$ for all $n,m$ and $i$, and
		\item $(p_{n;a}^{m;b})_{i,j}=(p_{a}^{b})_{i,j}$ for all $n,m$.
	\end{itemize}
	The above construction will then map such a 2-D QBD, considered for example in Ozawa~\cite{20190}, to a 1-D LD-QBD $\{(Z(t),\chi(t)):t\geq 0\}$ with parameters ${\bf P}_{n;a}^{m;b}={\bf P}_{a}^{b}$ and $\bLambda^{(n,m)}=\bLambda=diag(\lambda_i)_{i\in\{1,\ldots,k\}}$ for all $n,m$. The advantage of the one-dimensional representation is that the block matrices ${\bf Q}^{[n,n^{'}]}=[{\bf Q}^{[n,n^{'}]}_{i,j}]_{i,j\in\{1,2,\ldots ,m\}}$ are of finite size, and so the analysis can be performed using the existing results in the literature for the LD-QBDs, as we will discuss in our paper.
\end{Remark}

Using this new representation we will analyze the transient and steady-state behaviour of the random walk  considered here, in our future work.

\section{Convergence to a 2-D SFM}

Below we construct a sequence of two-dimensional QBDs $\left(\{(X_1^{(k)}(t),X_2^{(k)}(t),\varphi(t)):t\geq 0\}\right)_{k=1,2,\ldots}$ which converges in the distribution to $\{(Y_1(t),Y_2(t),\varphi(t)):t\geq 0\}$, a process which is a two-dimensional stochastic fluid model (2-D SFM).

This construction is a generalisation of the spatially coherent uniformisation of a SFM to a QBD proposed in~\cite{2013BO_uni}, to a 2-D level case studied here. The results follow directly by arguments in~\cite{2013BO_uni}. That is, the construction described below is a spatially coherent uniformisation of a 2D-SFM to a 2D-QBD (which can then be analysied using the methods of the earlier sections).

Let $\{\varphi(t):t\geq 0\}$ be a continuous-time Markov chain with state space $\mathcal{S}=\{1,2,\ldots,k\}$ and generator ${\bf T}=[\mathcal{T}_{ij}]_{i,j\in\mathcal{S}}$. Suppose that $c_{1,i}$ and $c_{2,i}$ is the collection of some real-valued rates determined for all $\varphi(t)=i\in\mathcal{S}$ and let $\mathcal{S}=\bigcup_{a,b\in\{+,-,0\}} \mathcal{S}_{a,b} $ be a partitioning of $\mathcal{S}$ such that $\mathcal{S}_{++}=\{i:c_{1,i}>0,c_{2,i}>0\}$, $\mathcal{S}_{+-}=\{i:c_{1,i}>0,c_{2,i}<0\}$, $\mathcal{S}_{+0}=\{i:c_{1,i}>0,c_{2,i}=0\}$, $\mathcal{S}_{-+}=\{i:c_{1,i}<0,c_{2,i}>0\}$, $\mathcal{S}_{--}=\{i:c_{1,i}<0,c_{2,i}<0\}$, $\mathcal{S}_{-0}=\{i:c_{1,i}<0,c_{2,i}=0\}$, $\mathcal{S}_{0+}=\{i:c_{1,i}=0,c_{2,i}>0\}$, $\mathcal{S}_{0-}=\{i:c_{1,i}=0,c_{2,i}<0\}$, and $\mathcal{S}_{00}=\{i:c_{1,i}=0,c_{2,i}=0\}$.

Then a 2-D SFM $\{(Y_1(t),Y_2(t),\varphi(t)):t\geq 0\}$ is a process with a two-dimensional level variable $(Y_1(t),Y_2(t))\in [0,\infty)\times [0,\infty)$ driven by the one-dimensional level variable $\varphi(t)$ which evolves according to the following assumptions: 
\begin{itemize}
	\item When $Y_1(t)>0$ then $dY_1(t)/dt=c_{1,\varphi(t)}$. When $Y_1(t)=0$ then $dY_1(t)/dt=max\{0,c_{1,\varphi(t)}\}$.
	\item When $Y_2(t)>0$ then $dY_2(t)/dt=c_{2,\varphi(t)}$. When $Y_2(t)=0$ then $dY_2(t)/dt=max\{0,c_{2,\varphi(t)}\}$.
\end{itemize}
That is, the process moves on the quadrant in a way such that the position $(Y_1(t),Y_2(t))$ is changing in a continuous manner depending on $\varphi(t)$.

We divide the quadrant into rectangles by discretising the $x$-axis with $\Delta x>0$ and the $y$-axis with $\Delta y>0$ such that $\Delta z=\sqrt{(\Delta x)^2+(\Delta y)^2}=1/k$ for some $k=1,2,\ldots$.

Consider $\left(\{(X_1^{(k)}(t),X_2^{(k)}(t),\varphi(t)):t\geq 0\}\right)_{k=1,2,\ldots}$ with a two-dimensional level variable $(X_1^{(k)}(t),X_2^{(k)}(t))$ taking values in $\{0,1,2,\ldots\}\times \{0,1,2,\ldots\}$ and one-dimensional phase variable $\varphi(t)$ driven by the above continuous-time Markov chain with state space $\mathcal{S}$ and generator ${\bf T}$. Assume that the generator
\begin{equation*}
{\bf Q}^{(k)} =
[q^{(k)}_{(n,m,i)(n^{'},m^{'},j)}]
\end{equation*}
of the $k$-th process in the sequence is such that, with $\Delta z=1/k$, the rates $q^{(k)}_{(n,m,i)(n^{'},m^{'},j)}$ are given by
\begin{eqnarray*}
\left\{
\begin{array}{ll}
\mathcal{T}_{i,j}&
j\not= i;(n,m)=(n^{'},m^{'})\\[2ex]
\frac{\sqrt{
		c_{1,i}^2+c_{2,i}^2 
	}}{\Delta z} 
	&
	j=i;(n,m)\not=(n^{'},m^{'})=([n+a]^+,[m+b]^+)
	\end{array}
	\right.\
\end{eqnarray*}
and $q_{(n,m,i)(n,m,i)}=-\sum_{(n,m,i)\not=(n^{'},m^{'},j)}q_{(n,m,i)(n^{'},m^{'},j)}$.

That is, at the times the 2D-QBD $\{(X_1^{(k)}(t),X_2^{(k)}(t),\varphi(t)):t\geq 0\}$ transitions from state $(n,m,i)$ to another state, it either transitions from phase $i$ to some $j$ at rate $\mathcal{T}_{i,j}$ without changing the level, or moves from level $(n,m)$ to $(n^{'},m^{'})=([n+a]^+,[m+b]^+)$ when $i\in\mathcal{S}_{a,b}$ at rate $\frac{\sqrt{
		c_{1,i}^2+c_{2,i}^2 
	}}{\Delta z}$ without changing the phase.

Consider the Poisson process $\{(N_{i,\Delta z}(t)):t\geq 0\}$ with rate $\vartheta_i(\Delta z)=\frac{\sqrt{
		c_{1,\varphi(t)}^2+c_{2,\varphi(t)}^2 
	}}{\Delta z} $ for every $i\in\mathcal{S}$. Let $L_i(t,\Delta z)$ be the random variable defined by $L_{i,\Delta z}(t)=N_{i,\Delta z}(t)\Delta z$. Then by standard properties of a Poisson Process, it follows that
	\begin{eqnarray*}
		\mathbb{E}(L_{i,\Delta z}(t))&=&
		\left(\sqrt{
			c_{1,\varphi(t)}^2+c_{2,\varphi(t)}^2 
		}\right) t,\\
		\lim_{\Delta z\to 0^{+}}Var(L_{i,\Delta z}(t))&=&
		0,
	\end{eqnarray*}
	and so, by the Pythagorean theorem, when the phase process spends $t$ units of time in phase $i$, then the change in level $X_1(\cdot)$ is approximately $|c_{1,i}| t$ and the change in level $X_2(\cdot)$ is approximately $|c_{2,i}| t$ when $\Delta z=\sqrt{(\Delta x)^2+(\Delta y)^2}$ is very small.

	Furthermore, $\left(\{(X_1^{(k)}(t),X_2^{(k)}(t),\varphi(t)):t\geq 0\}\right)_{k=1,2,\ldots}$ converges in distribution to $\{(Y_1(t),Y_2(t),\varphi(t)):t\geq 0\}$, which we write as
	$$\{ (X_{1}^{(k)}(\cdot),X_{2}^{(k)}(\cdot), \varphi_k(\cdot))\} \Rightarrow (Y_1(\cdot),Y_2(\cdot),\varphi(\cdot))$$ in $D_{\subR\times\subR\times S}[0,\infty)$.\\

In the future paper, we will discuss the methodology for the evaluation of the various stationary and transient quantities of the random walk considered here, and illustrate the theory through numerical examples of application.

\bibliographystyle{abbrv}
\bibliography{refs_rw}

\end{document}